\documentstyle{amsppt}
\magnification=\magstep1
\overfullrule=0pt
\loadbold
	\def\dist{\mathop{\text{\rm dist}}\nolimits} 
	\def\Int{\mathop{\text{\rm Int}}\nolimits} 
	\def\olim{\mathop{\overline{\text{\rm lim}}}\nolimits}
	\def\ulim{\mathop{\underline{\text{\rm lim}}}\nolimits}
	\def\oosc{\mathop{\overline{\text{\rm osc}}}\nolimits}
	\def\uosc{\mathop{\underline{\text{\rm osc}}}\nolimits}
	\def\os{\mathop{\text{\rm os}}\nolimits} 
	\def\osc{\mathop{\text{\rm osc}}\nolimits} 
	\def\sgn{\mathop{\text{\rm sgn}}\nolimits} 
	\def\supp{\mathop{\text{\rm supp}}\nolimits} 
	\def\varep{\varepsilon}
	\def\real{{\Bbb R}}
	\def\D{{\Cal D}}
	
	\def\G{{\Cal G}}
	\def\iitem{\itemitem}
	\def\chix{{\raise.5ex\hbox{$\chi$}}}
	\def\To{\Rightarrow}
	\def\btheta{{\boldsymbol\theta}}

\topmatter
\title On Functions of Finite Baire Index\endtitle
\author F. Chaatit*, V. Mascioni** and H. Rosenthal**\endauthor
\leftheadtext{F. Chaatit, V. Mascioni and H. Rosenthal}
\thanks * Some of the results given here forms part of the first named 
author's Ph.D.\ thesis at The University of Texas at Austin, supervised by 
the third named author.\endthanks
\thanks ** Research partially supported by NSF DMS-8903197.\endthanks 
\abstract 
It is proved that every function of finite Baire index on a separable metric 
space $K$ is a $D$-function, i.e., a difference of bounded semi-continuous 
functions on $K$. In fact it is a strong $D$-function, meaning it can be 
approximated arbitrarily closely in $D$-norm, by simple $D$-functions.
It is shown that if the $n^{th}$ derived set of $K$ is non-empty 
for all finite $n$, there exist 
$D$-functions on $K$ which are not strong $D$-functions. Further structural 
results for the classes of finite index functions and strong $D$-functions 
are also given.
\endabstract 
\endtopmatter

\document
\baselineskip=18pt

\head 1. Introduction\endhead 

Throughout, let $K$ be a separable metric space. 
A function $f:K\to \real$ is called a difference of bounded semi-continuous 
functions if there exist bounded lower semi-continuous functions $u$ and $v$ 
on $K$ with $f=u-v$. 
We denote the class of all such functions by $DBSC(K)$.  
We shall also refer to members of $DBSC(K)$ as {\it $D$-functions\/}. 
A classical theorem of Baire (cf.\ \cite{H, p.274}) 
yields that $f\in DBSC(K)$ if and only if 
there exists a sequence $(\varphi_j)$ of continuous functions on $K$ so that 
$$\sup_{k\in K} \sum |\varphi_j (k)| <\infty\quad\text{and}\quad 
f= \sum \varphi_j\ \text{point-wise.} 
\leqno(1)$$ 
Now defining $\|f\|_D = \inf \{\sup_{k\in K}\sum |\varphi_j|(k) :
(\varphi_j)$ is a 
sequence of continuous functions on $K$ satisfying (1)$\}$, 
it easily follows that $DBSC(K)$ is a Banach algebra; and of course 
$DBSC(K)\subset B_1(K)$ where $B_1(K)$ denotes the (bounded) first Baire 
class of functions  on $K$; i.e., the space of all bounded functions on $K$ 
which are the limit of a  point-wise convergent 
sequence of continuous functions on $K$. 

$DBSC(K)$ appears as a natural object in functional analysis. For example, if 
$X$ is a separable Banach space and $K$ is the unit ball of $X^*$ in the 
weak*-topology, then $X$ contains a subspace isomorphic to $c_0$ if and only 
if there is an $f$ in $X^{**}\sim X$ with $f\mid K$ in $DBSC(K)$ 
(cf.\ \cite{HOR}, \cite{R1}). 
Natural invariants for $DBSC(K)$ are used in a fundamental way in \cite{R1}, 
to prove that $c_0$ embeds in $X$ provided $X$ is non-reflexive and $Y^*$ is 
weakly sequentially complete for all subspaces $Y$ of $X$. 

We investigate here a special subclass of $DBSC(K)$, which we term 
$SD(K)$, and show that all functions of finite Baire index belong to 
this class. 

To motivate the definitions of these objects we first recall the following 
class of functions. 
Define $B_{1/2}(K)$ to be the set of all uniform limits of functions in 
$DBSC(K)$. 
(The terminology follows that in \cite{HOR}.) 
Functions in $B_{1/2}(K)$ may be characterized in terms of an intrinsic 
oscillation behavior, which we now give. 

For $f:K\to\real$ a given bounded function, let $Uf$ denote the upper 
semi-continuous envelope of $f$; $Uf(x) = \olim_{y\to x} f(y)$ for all 
$x\in K$.  (We use non-exclusive lim sups; thus equivalently, $Uf(x) = 
\inf_U \sup_{y\in U} f(y)$, the inf over all open neighborhoods of $x$.) 
Now we define $\uosc f$, the lower oscillation of $f$, by 
$$\uosc f(x) = \olim_{y\to x} |f(y) - f(x)|\ \text{ for all }\ x\in K\ .
\leqno(2)$$ 
Finally, we define $\osc f$, the oscillation of $f$, by 
$$\osc f = U \uosc f\ .
\leqno(3)$$ 

Now let $\varep>0$. We define the (finite) oscillation sets of $f$, 
$\os_j(f,\varep)$, as follows. 
Set $\os_0(f,\varep)= K$. 
Suppose $j\ge0$ and $\os_j(f,\varep)$ has been defined. 
Let $\os_{j+1}(f,\varep) = \{x\in L :\osc f\mid L(x)\ge \varep\}$, 
where $L= \os_j(f,\varep)$. 

We recall the following fact (\cite{HOR}). 

\proclaim{Proposition 1.1} 
Let $f:K\to\real$ be a given function. The following are equivalent: 
\vskip1pt 
\iitem{\rm1.} $f\in B_{1/2}(K)$.  
\iitem{\rm2.} For all $\varep>0$, there is an $n$ with $\os_n(f,\varep) 
= \emptyset$. 
\endproclaim 

\noindent (The proof given in \cite{HOR} for compact metric spaces works 
for arbitrary separable ones; cf.\ also \cite{R2}.) 

\demo{Remark} 
Actually, the sets defined in \cite{HOR} use what we term here the upper 
oscillation of $f$, defined by $\oosc f(x) = \olim_{y,z\to x} |f(y)-f(z)|$. 
It is easily seen that $\oosc f$ is upper semi-continuous and 
$$\tfrac12 \oosc f \le \osc f \le \oosc f\ .
\leqno(4)$$ 

Now define $K_j(f,\varep)$ inductively by 
$$K_0 (f,\varep) = K\quad\text{and}\quad 
K_{j+1}(f,\varep) = \{ x\in K_j :\oosc f\mid K_j (x)\ge \varep\}\ .$$ 
We then have by (4) that 
$$K_j (f,2\varep) \subset \os_j(f,\varep) \subset K_j(f,\varep) 
\ \text{ for all }\ j\ .
\leqno(5)$$ 
Thus $f$ satisfies 2 of 1.1 if and only if for all $\varep>0$, there is an 
$n$ with $K_n(f,\varep) =\emptyset$. 
\enddemo 

Proposition 1.1 suggests the following quantitative notion. 

\proclaim{Definition 1} 
Let $f:K\to\real$ be a given bounded function and $\varep>0$. 
We define $i(f,\varep)$, 
the $\varep$-oscillation index of $f$, to be $\sup \{n:\os_n(f,\varep) 
\ne\emptyset\}$. 
\endproclaim 

Thus Proposition 1.1 says that $f\in B_{1/2}(K)$ if and only if $i(f,\varep) 
<\infty$ for all $\varep>0$. 

\proclaim{Definition 2} 
A bounded function $f:K\to\real$ is said to be of finite Baire index if 
there is an $n$ with $\os_n(f,\varep) =\emptyset$ for all $\varep>0$. 
We then define $i(f)$, the oscillation index of $f$, by 
$$i(f) = \max_{\varep>0} i(f,\varep)\ .$$ 
\endproclaim 

\noindent Evidently $f$ is continuous if and only if $i(f)=0$. 

\demo{Remark} 
In \cite{HOR}, an index $\beta (f)$ is defined as $\beta (f)= \sup_{\varep>0} 
\min \{j:K_j(f,\varep) = \emptyset\}$. 
It follows from the remark following Proposition~1.1 that $f$ is of finite 
index if and only if $\beta (f)<\infty$, and then in fact $\beta (f) = i(f)+1$. 
\enddemo 

In \cite{HOR}, it is proved that finite index functions belong to 
$B_{1/4}(K)$, a class properly containing the $D$-functions. 
We obtain here that every function of finite Baire index belongs to $DBSC(K)$.  
In fact, we show that it belongs to the following subclass: 

\proclaim{Definition 3} 
A function $f:K\to\real$ is said to be a strong $D$-function if there exists 
a sequence $(\varphi_n)$ of simple $D$-functions with $\|f-\varphi_n\|_D\to0$.  
We denote the class of all strong $D$-functions by $SD(K)$. 
\endproclaim 

We may thus formulate one of our main results as follows:  

\proclaim{Theorem 1.2} 
Let $f:K\to\real$ be a function of finite Baire index. 
Then $f$ belongs to $SD(K)$. 
\endproclaim 

As we show below it is easily seen that every simple $D$-function has finite 
Baire index. Thus Theorem~1.2 yields that $SD(K)$ equals the closure, in 
$D$-norm, of the functions of finite index on $K$. 
Our proof essentially proceeds from first principles. 
An alternate argument, using transfinite oscillations, is given in \cite{R2}. 

An interesting special case of 1.2: 
{\it Let $f:[0,1]\to\real$ be bounded such that $\lim_{y\uparrow x} f(y)$, 
$\lim_{y\downarrow x} f(y)$ exist for all $x$. 
Then $f$ is in $SD[0,1]$.} 
The fact that such functions are in $DBSC[0,1]$ was initially proved jointly 
by the first and third named authors, and precedes the work given here 
\cite{C}. (It is a standard elementary result that if $f$ has these 
properties, then $\os_1(f,\varep)$ is {\it finite\/} for all $\varep>0$, 
hence $i(f)=1$.) 

It is evident that the simple $D$-functions form an algebra, hence 
$SD(K)$ is a Banach algebra. It is proved in \cite{R2} that $SD(K)$ is a 
lattice, i.e., $|f|\in SD(K)$ if $f\in SD(K)$. 
We prove here that the functions of finite index form an algebra and a 
lattice. This follows immediately from the following result. 

\proclaim{Theorem 1.3} 
Let $f,g$ be bounded real-valued functions on $K$, of finite index. 
Let $h$ be any of the functions $f+g$, $f\cdot g$, $\max \{f,g\}$, 
$\min \{f,g\}$. Then 
$$i(h) \le i(f) + i(g)\ . 
\leqno(6)$$ 
\endproclaim 

It is evident that if $f$ is of finite index, then for any non-zero scalar 
$\lambda$, $i(\lambda f) = i(f)$; also it is easy to show that 
$i(|f|) \le i(f)$. 
However the assertions of Theorem~1.3 appear to lie below the surface. 
The quantitative result which does the job (Theorem~2.8 below), is then 
applied to 
yield a necessary condition for a function to be in $SD(K)$, 
which is also sufficient in the case of upper semi-continuous functions. 

\proclaim{Theorem 1.4} 
Let $f:K\to\real$ be a given bounded function.  
\vskip1pt
(a) If $f\in SD(K)$, then 
$$\lim_{\varep\to0} \varep i(f,\varep) = 0 
\leqno(7)$$
\indent (b) If $f$ is semi-continuous and satisfies $(7)$, then 
$f\in SD(K)$. 
\endproclaim 

It is proved in \cite{R2} that every $SD$-function is a difference of 
strong $D$-semi-continuous functions. 
Evidently Theorem~1.4 yields an effective criterion for distinguishing 
the class of strong-$D$ semi-continuous functions. However 
one may construct functions, e.g., on $K= \omega^\omega+1$, which are 
not $D$-functions but satisfy (7), or which are $D$-functions but not 
$SD$-functions, and still satisfy (7).
An effective intrinsic criterion involving the ``$\omega^{th}$ oscillation'',  
which does distinguish $SD$-functions from $D$-functions, is given 
in \cite{R2}. 

We conclude the article by applying Theorem 1.4(a) to show that 
$DBSC(K)\sim SD(K)$ is non-empty for all interesting $K$. 

\proclaim{Proposition 1.5} 
Assume that $K^{(j)}$, the $j^{th}$ derived set of $K$, is non-empty for 
all $j=1,2,\ldots$. 
There exists a function $f$ on $K$ which is in $DBSC(K)$ but not in 
$SD(K)$. 
\endproclaim 

(An alternate proof of 1.5, using transfinite oscillations, is given in 
\cite{R2}.)

Recall that $K^{(j)}$ is defined inductively: 
For $M$ a topological Hausdorff space, let $M'$ denote the set of cluster 
points of $M$. Let $K^{(0)} = K$ and $K^{(j+1)} = (K^{(j)})'$ for all $j$. 
Now if $K$ fails the hypotheses of 1.5 there is an integer $n$ with 
$K^{(n+1)} = \emptyset$. 
Then every bounded function on $K$ is of index at most $n$, 
hence belongs to $SD(K)$. 
It can also be shown that if $K$ satisfies the hypotheses of 1.5, there 
exists an $f\in B_{1/2}(K)\sim DBSC(K)$, and also an $f\in B_1(K) \sim 
B_{1/2}(K)$.

\head Section 2.\endhead 

We begin with some preliminary results. 

\proclaim{Lemma 2.1} 
Let $f$ be a bounded non-negative lower semi-continuous function on $K$. 
Then $f\in DBSC(K)$ and $\|f\|_D = \|f\|_\infty$. 
Hence if $f$ is bounded semi-continuous, $\|f\|_D \le 3\|f\|_\infty$. 
\endproclaim 

\demo{Proof} 
By a classical result of Baire (cf.\ \cite{H}), there exists a sequence 
$(\varphi_j)$  of continuous functions on $K$ with $0\le \varphi_1\le 
\varphi_2\le \cdots$ and $\varphi_j \to f$ pointwise. 
Setting $u_1=\varphi_1$, $u_j = \varphi_j - \varphi_{j-1}$ for $j>1$, 
we have that $u_j\ge0$ for all $j$ and $\sum u_j=f$ point-wise. 
Thus $\|f\|_D \ge \|f\|_\infty$; the reverse inequality is trivial.

To see the last statement, let e.g., $f$ be bounded upper semi-continuous, 
$\lambda = \|f\|_\infty$, and note that $\lambda-f$ is non-negative lower 
semi-continuous. 
Thus $\|\lambda-f\|_D = \|\lambda-f\|_\infty \le 2\lambda$, 
so $\|f\|_D \le \lambda + \| \lambda-f\|_D\le 3\lambda$.\qed 
\enddemo 

\demo{Remark} 
It thus follows that if $f$ is a $D$-function, then $\|f\|_D = \inf 
\{ \|u+v\|_\infty :u,v\ge0$ are bounded lower semi-continuous with 
$f= u-v\}$. 
\enddemo 

Of course it follows immediately from Lemma 2.1 that if $U$ is an open 
non-empty subset of $K$, then  $\|\chix_U\|_D=1$, for $\chix_U$ is lower 
semi-continuous. 
In this case, the sequence $(\varphi_j)$ mentioned above can be easily chosen, 
using Urysohn's lemma. 
Indeed, if $U$ is closed, this is trivial. 
Otherwise, let $\varep_0>0$ be such that $\dist (x_0,\partial U) >\varep_0$ 
for some $x_0\in U$; set $F_n = \{x\in U :\dist (x,\partial U) \ge {\varep_0
\over n}\}$. 
Then $U= \bigcup_{j=1}^\infty F_j$ and for all $j$, $F_j$ is closed, 
$F_j \subset \Int F_{j+1}$. 
Now choose $[0,1]$-valued continuous functions $(\varphi_j)$ on $K$ so that 
for all $j$, $\varphi_j=1$ on $F_j$ and $\overline{\{x:\varphi_j(x)\ne 
0\}} \subset \Int F_{j+1}$. 
Then $\varphi_j \to \chix_U$ pointwise. 

Evidently it follows that if $W$ is a closed subset of $K$, then $\|\chix_W
\|_D \le 2$. In fact, if $W$ is a difference of closed sets; i.e., 
$W= W_1\sim W_2$, with $W_i$ closed for $i=1,2$, we again have that 
$\|\chix_W\|_D\le 2$, for $\|\chix_W\|_D\le \|\chix_{W_1}\|_D 
\|\chix_{\sim W_2}\|_D \le 2\cdot 1 = 2$. 

The following result shows that the simple $D$-functions are precisely 
those functions built up from the differences of closed sets. 

\proclaim{Proposition 2.2} 
Let $f$ be a simple real-valued function on $K$. 
The following are equivalent: 
\vskip1pt
\iitem{\rm 1)} $f\in B_{1/2}(K)$;  
\iitem{\rm 2)} $f$ is of finite Baire index;  
\iitem{\rm 3)} $f\in DBSC(K)$;  
\iitem{\rm 4)} There exist disjoint differences of closed sets $W_1,\ldots,W_m$ 
and scalars $c_1,\ldots,c_m$ with 
$$f= \sum_{i=1}^m c_i \chix_{W_i}\ .$$
\endproclaim 

\demo{Proof} 
Let us suppose $f$ is non constant, let $r_1,\ldots r_k$ be the distinct 
values of $f$, and set $\varep = \min \{ |r_i-r_j| : i\ne j$, 
$1\le i,j\le k\}$.  
Now if $W$ is a non-empty subset of $K$, $w\in W$, and $\osc  f\mid W(w) < 
\varep$, then $f\mid W$ is continuous at $w$; in fact there is an open 
neighborhood $U$ of $w$ with $f(x) = f(w)$ for all $x\in U\cap W$. 

Now suppose 1) holds, and let $n= i(f,\varep)$. 
By Proposition~1.1, $n<\infty$. 
We then obtain that defining $K_0 = K$ and $K_{j+1}=\{x\in K_j :f|K_j$ 
is discontinuous at $x$\}, for $1\le j\le n+1$, then $K_{n+1}=\emptyset$ and 
if $0 < \varep' \le\varep$, $\os_j (f,\varep') = K_j$ 
for all $1\le j\le n$. 
Hence in fact $i(f) = i(f,\varep) = n$, so 2) is proved. 
Of course 2) implies 1) by Proposition~1.1. 

It remains only to show that $1)\To 4)$, for evidently $4) \To 3) \To 1)$. 
Now fixing $0\le j\le n$, we have that 
$f$ is continuous on $K_j\sim K_{j+1}$. 
Let then $\ell = \ell (j)$ and 
$r_1^j,\ldots, r_\ell^j$ be the distinct values of $f$ on $K_j\sim K_{j+1}$; 
let $W_i^j = \{x\in K_j \sim K_{j+1}: f(x) = r_i^j\}$. 
Then $W_i^j$ is a clopen subset of $K_j\sim K_{j+1}$; it follows easily that 
in fact $W_i^j$ is then again a difference of closed sets in $K$, for all $i$, 
$1\le i\le \ell$, and thus 
$$ f= \sum_{j=0}^n \sum_{i=1}^{\ell(j)} r_i^j \chix_{W_i^j}\ , $$ 
proving 4).\qed
\enddemo 

\demo{Remark} 
The above proof yields that moreover if $W\subset K$, and $\chix_W$ is a 
$D$-function, then $W$ is a (disjoint) finite union of differences of closed 
sets; the converse is again immediate. 
This condition is incidentally equivalent to the condition that $W$ belongs 
to the algebra $\D$ of sets generated by the closed subsets of $K$. 
\enddemo 

We give some more preliminary results, before passing to the proof of 
Theorem~1.2. For $f:K\to \real$, we set $\supp f= \{k\in K: f(k)\ne0\}$. 
If $W\subset K$, we say that $f$ is supported on $W$ if $\supp f\subset W$. 

\proclaim{Lemma 2.3} 
Let $U$ be a non-empty open subset of $K$, and $f$ a bounded function on $K$, 
supported and continuous on $U$. Then $f\in SD(K)$ and $\|f\|_D=\|f\|_\infty$. 
\endproclaim 

\demo{Proof} 
Let us first show the norm identity. 
Note that since $f$ is bounded, if $u$ is a continuous function on $K$ with 
$u(x) = 0$ for all $x\notin U$, then $f\cdot u$ is continuous on $K$. 
Now choose $u_1,u_2,\ldots$ continuous non-negative functions on $K$ with 
$\chix_U = \sum u_j$ point-wise. 
But then $f= \sum f\cdot u_j$ point-wise, $f\cdot u_j$ is continuous on $K$ 
for all $j$, and $\sum |fu_j| \le \|f\|_\infty \sum u_j \le \|f\|_\infty$, 
so $\|f\|_D \le \|\sum |fu_j|\,\|_\infty \le \|f\|_\infty$; the reverse 
inequality is trivial. 

To see that $f$ is a strong  $D$-function, assume without loss of generality 
that $\|f\|_\infty=1$. Now fix $n$ a positive integer, and for each $j$, 
$-n\le j\le n$, define $K_j^n$ by 
$$K_j^n = \Big\{ x\in U : {j\over n} \le f(x) < {j+1\over n}\Big\}\ .
\leqno(8)$$ 
Finally, define $\varphi_n$ by 
$$\varphi_n = \sum_{j=-n}^n {j\over n} \chix_{K_j^n}\ . 
\leqno(9)$$ 
Then evidently by  the continuity of $f$, $K_j^n$ is a difference of closed 
sets in $U$, and hence in $K$, for all $j$, so $\varphi_n$ is a simple 
$D$-function; moreover we have 
$$0\le f-\varphi_n \le {1\over n}\ . 
\leqno(10)$$ 

Thus to  show that $\|f-\varphi_n\|_D \to0$ as $n\to\infty$, we need only show 
that $f-\varphi_n$ is lower semi-continuous; for then 
$\|f-\varphi_n\|_D \le {1\over n}$ by (10) and  Lemma~2.1. 

Let $\psi = f-\varphi_n$, and suppose it were false that $\psi$ is lower 
semi-continuous. 
We may then choose $x\in K$ and $(x_m)$ a sequence in $K$ with $x_m\to x$ 
so that $(\psi (x_m))$ converges and 
$$\lim_{m\to\infty} \psi (x_m) < \psi (x)\ . 
\leqno(11)$$ 

Evidently then $x\in U$, since $x\notin U$ implies $\psi (x)=0\le \psi (x_m)$ 
for all $m$. By passing to a subsequence, we may then assume without loss of 
generality that there is a $j$,  $-n\le j\le n$, with $x_m \in K_j^n$ 
for all $m$. 
But since  $f$ is continuous on $U$, $\lim_{m\to\infty} f(x_m) = f(x)$; 
if also $x\in K_j^n$, then since $\psi (x_m) = f(x_m) - {j\over n}$ for 
all $m$, we have that $\lim_{m\to\infty} \psi (x_n) = f(x)-{j\over n} = 
\psi (x)$, a contradiction. 
If $x\notin K_j^n$, by continuity of $f$ 
we must have that $f(x)={j+1\over n}$. 
But then $x\in K_n^{j+1}$, so $\psi (x) = 0<  {j+1\over n} - {j\over n} 
= \lim_{m\to\infty} \psi (x_m)$ again contradicting (11).\qed 
\enddemo 

Our next preliminary result deals with extension issues. 
(For $W\subset K$ and \hfill\break 
$f:W\to\real$, $f\cdot \chix_W$ denotes the function 
which is zero off $W$ and agrees with $f$ on $W$.) 

\proclaim{Lemma 2.4} 
Let $W\subset K$ be a difference of closed sets and $f$ in $DBSC(W)$. 
Then $f\cdot \chix_W$ is in $DBSC(K)$ and 
$$\| f\cdot \chix_W\|_{D(K)} \le 2\|f\|_{D(W)}\ ; 
\leqno(12)$$ 
if $W$ is an open set, then 
$$\|f\cdot\chix_W\|_{D(K)} = \|f\|_{D(W)}\ . 
\leqno(13)$$ 
Moreover if $f\in SD(W)$, then $f\chix_W \in SD(K)$.
\endproclaim 

\demo{Proof} 
Suppose first that $W$ is open, and let $(\varphi_j)$ in $C(K)$ be such that 
the $\varphi_j$'s are non-negative and $\sum \varphi_j = \chix_W$ 
point-wise. Let $\varep>0$ and choose $(\psi_j)$ in $C(W)$ with 
$\sum |\psi_j| < \|f\|_{D(W)} + \varep$ and $f= \sum \psi_j$ point-wise 
on $W$. Now identifying $\psi_j$ with $\psi_j\cdot \chix_W$, 
$\psi_j\cdot\varphi_i$ is continuous on $K$ for all $i$ and 
$j$, and we have that $\sum_{i,j} |\psi_j\varphi_i| = 
\sum_j |\psi_j|\chix_W \le \|f\|_{D(W)} +\varep$, with 
$\sum_{i,j}\psi_j\varphi_i = f\chix_W$. Thus  $\|f\chix_W\|_{D(K)} \le 
\|f\|_{D(W)} +\varep$ for all $\varep>0$; so $\|f\chix_W\|_{D(K)}\le 
\|f\|_{D(W)}$. 
The reverse inequality is trivial, so (13) is established. 

Next, suppose that $W$ is closed, and again let $\varep>0$. 
As noted following Lemma~2.1, we may choose $u,v$ non-negative lower 
semi-continuous on $W$ with 
$$f= u-v\quad\text{and}\quad \|u+v\|_\infty <\|f\|_{D(W)} +\varep\ . 
\leqno(14)$$
Now let $\lambda= \|u+v\|_\infty$ and let $\tilde u = \lambda\chix_{\sim W} + 
u\chix_W$, $\tilde v= \lambda \chix_{\sim W} + v\chix_W$. 
It follows easily that $\tilde u$ and $\tilde v$ are both non-negative 
lower semi-continuous on $K$ and of course 
$$f\chix_W = \tilde u - \tilde v\ ,\qquad \|\tilde u+\tilde v\|_\infty 
= 2\lambda\ .
\leqno(15)$$ 
Thus by the observation following Lemma 2.1, $\|f\cdot \chix_W\|_D \le 
2\lambda < 2\|f\|_{D(W)} + 2\varep$. 
Since $\varep>0$ is arbitrary, (12) is proved for closed $W$. 

Now suppose $W$ is a difference of closed sets. 
Choose $U$ open, $L$ closed with $W= U\cap L$. 
Then $W$ is a relatively closed subset of $U$, so we have that $f\cdot \chix_L
\mid U$ belongs to $DBSC(U)$ with $\|f\cdot \chix_L\mid U\|_{D(U)} \le 
2\|f\|_{D(W)}$. 
But then by (13), $f\cdot\chix_W =(f\cdot \chix_L)\mid U\cdot\chix_U$ belongs 
to $DBSC(K)$ and $\|f\cdot \chix_W\| \le \|f \cdot \chix_L\mid U\|_{D(W)} 
\le 2\|f\|_{D(W)}$, proving (12). 

Finally, suppose $f\in SD(W)$. 
Then given $\varep>0$, choose $g$ a simple $D$-function on $W$ with 
$$\|g-f\|_{D(W)} <\varep\ .
\leqno(16)$$ 

By Proposition 2.2, there are disjoint differences of closed sets in $W$, 
$W_1,\ldots,W_k$, and scalars $c_1,\ldots,c_k$ with $g= \sum_{i=1}^k c_i 
\chix_{W_i}$ on $W$. 
But then for all $i$, 
$W_i$ is actually a difference of closed sets in $K$, and thus 
$g\cdot\chix_W$ is a simple $D$-function on $K$. 
Then by (12), 
$$\|(g-f)\chix_W \| = \|g\chix_W - f\chix_W\| <2\varep\ . 
\leqno(17)$$ 
Thus the final assertion of the Lemma is established.\qed
\enddemo 

\demo{Remark} 
Using the comment following Proposition 2.2, 
we obtain that if $W\subset K$ is in 
$\D$ (i.e., $\chix_W$ is a $D$-function), then for $f:W\to \real$ a bounded 
function, $f$ is a $D$-function on $W$ if and only if $f\chix_W$ is a 
$D$-function on $K$; moreover $f\in SD(W)$ if and only if $f\chix_W \in 
SD(K)$. 
\enddemo 

Before giving the proof of Theorem 1.2, we recall the following standard 
result. 

\proclaim{Lemma 2.5} 
Let $\varep>0$, and suppose $f:K\to\real$ is such that $\osc f\le \varep$ 
on $K$. There exists $\varphi:K\to\real$ continuous with $|f-\varphi| 
\le\varep$ on $K$. 
\endproclaim 

\demo{Proof} 
Let $Lf$ be the lower semi-continuous envelope of $f$; 
$Lf(x) = \ulim_{y\to x} f(y)$ for all $x\in X$. 
Then we have that 
$$\oosc f = Uf - Lf\ .
\leqno(18)$$ 
Since $\oosc f\le 2\osc f$, $\oosc f\le 2\varep$ on $K$. 
Thus we have by assumption that 
$$Uf-\varep \le Lf +\varep\ . 
\leqno(19)$$ 
By the Hahn  interposition theorem (cf.\ \cite{H}, p.276), there exists 
$\varphi$ continuous with 
$$Uf-\varep \le \varphi \le Lf +\varep\ . 
\leqno(20)$$ 
Since $f\le Uf$ and $Lf \le f$, $\varphi$ satisfies the conclusion of 
the Lemma.\qed 
\enddemo 

We now treat the proof  of Theorem 1.2. 
It is convenient to consider a larger class; for $n\ge0$, let $\G_n$ 
denote the family of all bounded functions $f:K\to\real$ so that there 
exists an open set $U$ with $f$ supported on $U$ and $i(f\mid U) \le n$. 
The following quantitative result yields Theorem~1.2 immediately. 

\proclaim{Theorem 2.6} 
Let $n\ge0$ and $f\in \G_n$. 
Then $f\in SD(K)$ and\hfill\break 
$\|f\|_D\le (2^{n+1}-1)\|f\|_\infty$. 
\endproclaim 

\demo{Remark} 
Of course it follows {\it a-posteriori\/} 
that if we prove the result just for functions 
$f$ of index $n$, then it holds immediately for functions in $\G_n$, by 
Lemma~2.4. The class $\G_n$ is needed for our proof, however.  We also note 
that the argument given in \cite{R2}, using transfinite oscillations, gives 
the optimal estimate: if $i(f)\le n$, then $\|f\|_D\le (2n+1) \|f\|_\infty$. 
\enddemo 

We prove 2.6 by induction on $n$. 
The case $n=0$ follows immediately from Lemma~2.3. 
Now let $n>0$ and suppose 2.6 proved for ``$n$''~$=n-1$. 

\proclaim{Lemma 2.7} 
Let $f\in \G_n$ and $\varep>0$. There exist functions $g$ and $h$ 
with $f=g+h$, $g\in \G_n$, $h\in SD(K)$, and 
$$\|h\|_D \le (2^{n+1}-1)\|f\|_\infty\ ,\qquad \|g\|_\infty \le\varep\ . 
\leqno(21)$$ 
\endproclaim 

\demo{Proof} 
Let $\lambda_j = 2^{j+1}-1$ for $j=0,1,2,\ldots$.  
Let $U$ be chosen with $f$ supported in $U$ and $i(f\mid U) \le n$. 
Let $W= \{x\in U:\osc f(x) \ge\varep\}$. 
It follows that  $W$ is a relatively closed subset of $U$ and 
$$i(f\mid W)\le n-1\ .
\leqno(22)$$ 
Thus by induction hypothesis and Lemma 2.4, 
$$f\cdot \chix_W \in SD(K)\quad\text{and}\quad 
\|f\cdot \chix_W\|_D \le 2\lambda_{n-1} \|f\|_\infty\ .
\leqno(23)$$ 

Now by Lemma 2.5, we may choose $\varphi :U\sim W\to\real$, $\varphi$ 
continuous on $U\sim W$, with 
$$\|\varphi\|_\infty \le \|f\|_\infty\quad\text{and}\quad 
|\varphi (x) - f(x)| \le \varep \quad\text{for all}\quad 
x\in U\sim W\ ,
\leqno(24)$$ 
Indeed, 2.5 gives $\tilde \varphi$ with $\tilde\varphi$ continuous and 
$|\tilde\varphi-f| \le \varep$ on $U\sim W$. 
But simply define $\varphi (x) = \tilde\varphi(x)$ if $|\tilde\varphi (x)| 
\le \|f\|_\infty$, and $\varphi (x) = \|f\|_\infty \sgn f(x)$ otherwise. 

Let $g$ and $h$ be defined by 
$$g= (f-\varphi) \chix_{U\sim W}\quad ,\quad 
h= f\cdot \chix_W + \varphi\cdot\chix_{U\sim W}\ .
\leqno(25)$$ 

Now evidently $\supp g\subset U\sim W$; since $\varphi$ is continuous 
on $U\sim W$, it follows that $i((f-\varphi)\mid U\sim W) \le i(f\mid U)\le n$; 
hence $g\in \G_n$, and by (24), $\|g\|_\infty \le\varep$. 

Evidently, $f=g+h$; finally, by (23) and Lemma~2.3, $h\in SD(K)$ and 
$$\|h\|_D\le (2\lambda_{n-1} + 1)\|f\|_\infty = \lambda_n \|f\|_\infty\ .
\eqno\qed$$
\enddemo

\demo{Proof of Theorem 2.6 for $n$} 
Fix $\varep>0$. 
We may choose by induction sequences $(h_j)$ and $(g_j)$ so that for all $j$, 
$$\leqalignno{
&f= h_1+\cdots +h_j + g_j &\text{\rm (26i)}\cr
&h_j\in SD(K)\quad ,\quad  g_j\in \G_n &\text{\rm (26ii)}\cr
&\|h_1\|_D \le\lambda_n \|f\|_\infty\quad ,\quad 
\|h_j\|_D \le {\varep\over 2^{j-1}} \ \text{ for }\ j>1 &\text{\rm (26iii)}\cr
&\|g_j\|_\infty \le {\varep\over \lambda_n 2^j}\ .&\text{\rm (26iv)}\cr}$$ 

Indeed, by Lemma 2.7, we may choose $h_1\in SD(K)$ and $g_1\in \G_n$ 
with $f= h_1+g_1$, $\|h_1\|\le\lambda_n\|f\|_\infty$, $\|g_1\|_\infty 
\le {\varep\over \lambda_n 2}$. 

Now suppose $j\ge1$ and $h_1,\ldots,h_j$, $g_j$ chosen satisfying 
(26i)--(26iv). Since $g_j\in \G_n$, by Lemma~2.7 we may choose $h_{j+1}\in 
SD(K)$ and $g_{j+1}\in \G_n$ with $g_j = h_{j+1}+ g_{j+1}$,  
$$\|h_{j+1}\|_D \le \lambda_n \|g_j\|_\infty\quad\text{and}\quad 
\|g_{j+1}\|_\infty \le {\varep\over\lambda_n 2^{j+1}}\ . 
\leqno(27)$$ 
Then (26i)--(26iv) hold at $j+1$. 

Since the $D$-norm is trivially larger than the sup-norm and $\|g_j\|_\infty 
\to0$, it follows from (26i) and (26iii) that $\sum h_i$ converges 
uniformly to $f$. Since $DBSC(K)$ is a Banach space, 
$\sum\|h_j\|_D< \infty$, and $h_j\in SD(K)$ for all $j$, 
it follows that $f\in SD(K)$. Finally, we have by (26iii) that 
$$\|f\|_D \le \lambda_n\|f\|_\infty + \sum_{j=2}^\infty {\varep\over 2^{j-1}} 
= \lambda_n \|f\|_\infty +\varep\ . 
\leqno(28)$$ 
Since $\varep>0$ is arbitrary, Theorem 2.6 is proved.\qed 
\enddemo 

We turn now to Theorem 1.3. This follows immediately from the following 
result. 

\proclaim{Theorem 2.8} 
Let $f,g\in B_{1/2}(K)$, and $\varep>0$. 
Then the following hold. 
\vskip1pt
\iitem{\rm (a)} $i(f+g,\varep) \le i(f,{\varep\over2}) + i(g,{\varep\over2})$.  
\iitem{\rm (b)} $i(f\cdot g,\varep) \le i(f,{\varep\over 2G}) 
+ i(g,{\varep\over 2F})$ where $F=\|f\|_\infty$, $G= \|g\|_\infty$, 
and it is assumed that $F,G>0$.
\iitem{\rm (c)} $i(h,\varep) \le i(f,\varep) + i(g,\varep)$ where 
$h= f\vee g$ or $h=f\wedge g$.
\endproclaim 

We give the detailed proof of (a) (which is also needed later), and then  
indicate how (b), (c) follow by the same method. 

We first note the following fact. 

\proclaim{Lemma 2.9} 
Let $W_1,\ldots,W_n$ be closed non-empty sets with 
$K= \bigcup_{i=1}^n W_i$ and $f:K\to \real$ a bounded function. Then 
$$\osc f = \max_{1\le i\le n} (\osc f\mid W_i) \chix_{W_i}\ . 
\leqno(29)$$ 
\endproclaim 

\demo{Proof} 
We first note that 
$$\uosc f = \max_{1\le i\le n} (\uosc f\mid W_i) \chix_{W_i}\ . 
\leqno(30)$$ 
For let $x\in K$ and choose $(x_m)$ in $K$ with $x_m\to x$ and 
$\uosc f(x) = \lim_{n\to\infty} |f(x_n)-f(x)|$. 
We may choose $i$ and $m_1<m_2<\cdots$ with $x_{m_j} \in W_i$ for all $j$. 
But then $x\in W_i$ and so $\uosc f(x) \le \uosc f\mid W_i(x) 
\le \max_\ell (\uosc f\mid W_\ell) \chix_{W_\ell} (x)$. 
The reverse inequality is trivial, so (30) follows. 

Now again let $x\in K$ and choose $(x_m)$ in $K$ with $x_m\to x$ and 
$\osc f(x) = \lim_{n\to\infty} \uosc f(x_m)$. 
By (30), we may again choose $m_1<m_2<\cdots$ and $i$ with 
$\uosc f(x_{m_j}) = \uosc f\mid W_i \chix_{W_i} (x_{m_j})$ for all $j$. 
Now if $\osc f(x)=0$, (29) is trivial. 
Otherwise, without loss of generality, $\osc f(x_{m_j}) > 0$ for all $j$; 
hence then $x_{m_j} \in W_i$ and so $x\in W_i$, whence 
$\osc f(x) \le \osc f\mid W_i(x) \le \max_\ell (\osc f\mid W_\ell) 
\chix_{W_\ell} (x)$. 
Again the reverse inequality is trivial, so (29) holds.\qed
\enddemo 

Now let $f,g$ be as in Theorem 2.8, and $\varep>0$ be given. 
For each $n=1,2,\ldots$ and ${\btheta} =  (\theta_1,\ldots,\theta_n)$ 
with $\theta_i = 0$ or $1$ for all $1\le i\le n$, we define closed subsets 
$L({\btheta})$ of $K$ as follows: 
$$L(0) = \Bigl\{ x\in K: \osc f(x) \ge {\varep\over2}\Bigr\}\quad ;\quad 
L(1) = \Bigl\{ x\in K :\osc g(x)\ge {\varep\over2}\Bigr\}\ .
\leqno(31)$$ 
If $n\ge 1$ and $L(\btheta) = L(\theta_1,\ldots,\theta_n)$ is defined, let 
$$\left\{ \eqalign{
&L(\theta_1,\ldots,\theta_{n+1}) = \Bigl\{ x\in L(\btheta) :\osc f\mid 
L(\btheta) \ge {\varep\over2}\Bigr\}\ \text{ if }\ \theta_{n+1}=0\cr 
&L(\theta_1,\ldots,\theta_{n+1}) = \Bigl\{ x\in L(\btheta) :\osc g\mid 
L(\btheta) \ge {\varep\over2}\Bigr\}\ \text{ if }\ \theta_{n+1}=1\ .\cr}
\right.\leqno(32)$$ 

These sets are closed, since $\osc f$, $\osc g$ are upper semi-continuous 
functions. We then have for all $n$ that 
$$\os_n (f+g,\varep) \subset \bigcup_{\btheta \in \{0,1\}^n} L(\btheta)\ . 
\leqno(33)$$ 

We prove this by induction on $n$. 
Now for $n=1$, since it is easily seen that $\osc(f+g) \le \osc f+ \osc g$, 
we then have that $\osc (f+g)(x)\ge\varep$ implies $\osc f(x)\ge {\varep 
\over 2}$ or $\osc g(x)\ge {\varep\over2}$; this gives 
$\os_1(f+g,\varep) \subset L(0)\cup L(1)$. 
Suppose (33) is proved for $n$, and suppose $K_n = \osc_n(f+g,\varep)$ and 
$x\in \os_{n+1}(f+g,\varep)$. 
Thus $\osc (f+g)\mid K_n(x)\ge\varep$. 
By the preceding lemma and (33), we may then choose $\btheta \in \{0,1\}^n$ 
with $x\in K_n\cap L(\btheta)$ and 
$$\eqalign{\osc (f+g)\mid K_n(x)  & = \osc (f+g)\mid K_n\cap L (\btheta)(x)\cr 
&\le \osc (f+g)\mid L(\btheta) (x)\cr 
&\le \osc f\mid L(\btheta)(x) + \osc g\mid L(\btheta)(x)\ .\cr}$$
It follows immediately that $x\in L(\theta_1,\ldots,\theta_n,0) \cup 
L(\theta_1,\ldots,\theta_n,1)$; thus (32) holds at $n+1$. 

Next, fix $n$ and $\btheta \in \{0,1\}^n$. Let 
$$j = j(\btheta)= \#\, \{1\le i\le n: \theta_i = 0\}\quad ,\quad 
k = k(\btheta) = \#\, \{1\le i\le n: \theta_i = 1\}\ . 
\leqno(34)$$ 
Then we claim  
$$L(\btheta) \subset \os_j\left( f,{\varep\over2}\right)\cap 
\os_k \left( g,{\varep\over2}\right)\ . 
\leqno(35)$$ 

Again we prove this by induction on $n$. 
The case $n=1$ is trivial, by the definitions of $L(0)$ and $L(1)$. 
Now suppose (35) is proved for $n$, and $(\theta_1,\ldots,\theta_{n+1})$ 
is given; let $j= j(\theta_1,\ldots,\theta_n)$ and $k=k(\theta_1,\ldots,
\theta_n)$. 
Now if $\theta_{n+1}=0$, then $j(\theta_1,\ldots,\theta_{n+1}) = j+1$ 
and $k(\theta_1,\ldots,\theta_{n+1})=k$; then by (35), 
$L(\theta_1,\ldots,\theta_{n+1}) \subset L(\theta_1,\ldots,\theta_n) 
\subset \os_k(g,{\varep\over2})$ and by definition and (35), 
$$\eqalign{L(\theta_1,\ldots,\theta_{n+1}) &\subset \left\{ x\in \os_j 
\Bigl( f,{\varep\over2}\Bigr) : \osc f\mid \os_j\Bigl( f,{\varep\over2}
\Bigr) (x) \ge {\varep\over2}\right\}\cr 
& = \os_{j+1}\Bigl( f,{\varep\over2}\Bigr) \ .\cr}$$ 
Of course if $\theta_{n+1}=1$, we obtain by the same reasoning that 
$L(\theta_1,\ldots,\theta_{n+1}) \subset \os_j(f,{\varep\over2})\cap 
\os_{k+1}(g,{\varep\over2})$ and $j=j(\theta_1,\ldots,\theta_{n+1})$, 
$k+1 = k(\theta_1,\ldots,\theta_{n+1})$; 
thus (35) is proved for $n+1$, and so established for all $n$ by 
induction. 

Now suppose, for a given $n$, that $\os_n(f+g,\varep) \ne\emptyset$. 
Then by (33), there is a $\btheta \in \{0,1\}^n$ with 
$L(\btheta) \ne\emptyset$. Thus letting $j$ and $k$ be as in (34), 
we have by (35) that $\os_j(f,{\varep\over2})\ne\emptyset$ and 
$\os_k(g,{\varep\over2}) \ne\emptyset$. 
But then $n=j+k\le i(f,{\varep\over2}) + i(g,{\varep\over2})$. 
Theorem~2.8(a) is thus established. 

To see 2.8(b), note for any $y$ and $x\in K$ that 
$$|f(y)g(y)-f(x)g(x)| \le G|f(y) - f(x) | + F|g(y)-g(x)|\ .
\leqno(36)$$ 
Hence we have that fixing $x\in K$, then 
$\uosc fg(x)\le G\uosc f(x) + F\uosc g(x)$, whence 
$$\osc fg(x) \le G\osc f(x) + F\osc g(x)\ . 
\leqno(37)$$ 
Thus $\osc fg(x)\ge\varep$ implies $\osc f(x)\ge {\varep\over 2G}$ or 
$\osc g(x)\ge {\varep\over2F}$. 
We now prove (b) by defining the sets $L(\btheta)$ by 
$L(0) = \os_1 (f,{\varep\over2G})$, 
$L(1) = \os_1 (g,{\varep\over2F})$, 
and for $\btheta = (\theta_1,\ldots,\theta_{n+1})$, 
$L(\theta_1,\ldots,\theta_{n+1}) = \{x\in L(\btheta): \osc f\mid L(\btheta) 
\ge {\varep\over2G}\}$ if $\theta_{n+1}=0$, and 
$L(\theta_1,\ldots,\theta_{n+1}) = \{x\in L(\btheta) :\osc g\mid L(\btheta) 
\ge {\varep\over2F}\}$ if $\theta_{n+1}=1$. 
Then we proceed exactly as in case (a). 
Finally, for case (c), we note that if $h$ is as in (c) and $x\in K$, then
$$\osc h(x)\ge\varep \ \text{ implies }\ \osc f(x) \ge\varep \text{ or } 
\osc g\ge \varep\ .
\leqno(38)$$ 

Suppose this were false. Then we can choose $0<\varep' <\varep$ and $U$ 
an open neighborhood of $x$ with 
$$\osc f(u) < \varep'\ \text{ and }\ \osc g(u) < \varep'\ \text{ for all }\ 
u\in U\ . 
\leqno(39)$$ 
Now fix $u\in U$; we can then choose $V$ an open neighborhood of $u$ with 
$V\subset U$ and 
$$|f(v)-f(u)| <\varep'\ \text{ and }\  |g(v)-g(u)| <\varep'\ \text{ for all }
\ v\in V\ .
\leqno(40)$$ 
Suppose e.g., $h=f\vee g$ and $v\in V$ with $(f\vee g)(v)= f(v)$, 
$(f\vee g)(u)=g(u)$. But then by (40) and the above, 
$$f(v) \ge g(v) > g(u)-\varep'\ \text{ so }\ f(v) - g(u) >-\varep'
\leqno(41)$$ 
and 
$$f(v) < f(u) +\varep' \le g(u) +\varep'\ \text{ so }\ f(v)-g(u)<\varep'\ .
\leqno(42)$$ 
It thus follows from (40)--(42) that 
$$|h(v)-h(u)| <\varep'\ .
\leqno(43)$$ 
If e.g., $f\vee g(v) = f(v)$ and $f\vee g(u) = f(u)$, (43) follows 
immediately from (40), so (43) holds for all $v\in V$. 
Thus we obtain $\uosc h(u)\le \varep'$; but since $u\in U$ is 
arbitrary, we also have $\osc h(x)\le\varep'$, a contradiction. 
The proof for $h=f\wedge g$ is the same. 

Evidently (38) yields that 
$\os_1(h,\varep) \subset \os_1(f,\varep) \cup \os_1(g,\varep)$; 
we then proceed as in case~(a), except that the sets $L(\theta_1,\ldots,
\theta_n)$ are defined by replacing ``$\varep$'' by ``$\varep\over2$'' 
in (31), (32).\qed
\medskip

We next treat Theorem 1.4. We first recall the following fact. 

\proclaim{Lemma 2.9} 
Let $f\in D(K)$. Then $\varep i (f,\varep) \le 4\|f\|_D$. 
\endproclaim 

This follows immediately from the definitions, the fact that 
$\os_j(f,\varep) \subset K_j(f,\varep)$ for all $j$, and Lemma~2.4 
of \cite{HOR}. 
(A direct proof of 2.9 is given in \cite{R2} yielding the refinement 
that $\varep i(f,\varep) \le \|f\|_D$.) 

\demo{Proof of Theorem 1.4}
Suppose first that $f\in SD(K)$, $\eta>0$, and choose $g$ a simple 
$D$-function with $\|f-g\|_D\le \eta$. It then follows by Lemma~2.9 that 
$$\varep i (f-g,\varep) \le 4 \eta\ \text{ for all }\ \varep>0\ . 
\leqno(44)$$ 

Now since $g$ is a simple $D$-function, $g$ has finite index (by 
Proposition~2.2); say $N= i(g)$. Then by Theorem~2.8(a) and (44), for any 
$\varep>0$, 
$$\eqalign{\varep i(f,\varep) & \le \varep i\Bigl(f-g,{\varep\over2}\Bigr) 
+\varep i\Bigl( g,{\varep\over2}\Bigr) \cr 
&\le 8\eta + \varep N\ .\cr}$$ 
Hence $\olim_{\varep\to0} \varep i(f,\varep) \le 8\eta$. 
Since $\eta>0$ is arbitrary, (7) is proved. 

Finally, to prove (b) of Theorem 1.4, 
suppose without loss of generality that $f$ 
is upper semi-continuous and satisfies (7), let $\eta>0$, and choose 
$0<\varep<\eta$ with 
$$\varep i(f,\varep) < \eta\ . 
\leqno(45)$$ 

Let then $n= i(f,\varep)$ and set $K^j = \os_j(f,\varep)$ for all $j$. 
Thus $K^n \ne \emptyset$, $K^{n+1} = \emptyset$, and for $0\le j\le n$, 
$\osc (f\mid K^j\sim K^{j+1})<\varep$. 
Thus for all $j$, we may choose by Lemma~2.5 a continuous function  
$\varphi_j$ on $K^j\sim K^{j+1}$ with 
$$|\varphi_j -f| \le \varep \ \text{ on }\ K^j\sim K^{j+1}\ . 
\leqno(46)$$ 

Now set $g=\sum_{j=0}^n \varphi_j \chix_{K^j\sim K^{j+1}}$. 
By Lemmas~2.3 and 2.4, $g\in SD(K)$. 
Now fixing $j$ and letting $W= K^j \sim K^{j+1}$, then $(f-g)\mid W$ is upper 
semi-continuous, hence by Lemma~2.1 and (46), 
$$\|f-g\|_{D(W)} \le 3\|f-g\|_\infty \le 3\varep\ .
\leqno(47)$$ 

Then by Lemma 2.4, 
$$\|(f-g)\chix_W \|_{D(K)} \le 6\varep\ . 
\leqno(48)$$
Hence 
$$\eqalign{\|f-g\|_D & = \sum_{j=0}^n \|(f-g)\chix_{K^j\sim K^{j+1}}\|_D\cr 
&\le \sum_{j=0}^n \|(f-g)\chix_{K^j\sim K^{j+1}}\|_D\cr
&\le 6n\varep + 6\varep\cr
&< 7\eta \ \text{ by (45).}\cr}$$
Since $\eta>0$ is arbitrary and $SD(K)$ is closed in $DBSC(K)$, we obtain 
that $f\in SD(K)$, thus completing the proof of 
Theorem~1.4.\qed
\enddemo 

\demo{Remark} 
Define $B_{1/2}^0(K)$ to be the family of all bounded functions $f:K\to\real$ 
which satisfy (7). Evidently we have (by the preceding result) that 
$SD(K)\subset B_{1/2}^0(K) \subset B_{1/2}(K)$. 
We have moreover that $B_{1/2}^0(K)$ is an algebra and a lattice, 
by Theorem~2.8. As noted in the introduction, it can be shown that there are 
non-$D$-functions in $B_{1/2}^0(K)$, and also 
$(DBSC(K) \sim SD(K)) \cap B_{1/2}^0 (K)\ne\emptyset$ (for suitable $K$). 
It can be seen that $B_{1/2}^0(K)$ is a complete linear topological space 
under the quasi-norm $\|f\|= \sup_{\varep>0} \varep i(f,\varep) 
+ \|f\|_\infty$. 
\enddemo 

We finally consider Proposition 1.5. The construction uses some preliminary 
results. 

\proclaim{Lemma 2.10} 
Let $n\ge1$ and $K= K_0\supset K_1\supset \cdots\supset K_n$ be closed 
non-empty sets with $K_i$ nowhere dense relative to $K_{i-1}$ for all 
$1\le i\le n$. Also let $K_{n+1}=\emptyset$. 
Let $E= \bigcup_{0\le i\le [n/2]} K_{2i} \sim K_{2i+1}$. Then 
$$i(\chix_E) = i(\chix_E,\varep) = n\ \text{ for all }\ 0<\varep\le1\ . 
\leqno(49)$$ 
Moreover $\|\chix_E\|_D \le n+1$. 
\endproclaim 

\demo{Proof} 
Fix $0<\varep\le1$. We prove by induction on $j$ that 
$$\os_j(\chix_E,\varep) = K_j\ \text{ for all }\ 0\le j\le n\ . 
\leqno(50)$$ 
Then  since $\chix_E$ is constant on $K_n$, $\os_{n+1}(\chix_E,\varep) 
= \emptyset$, yielding (49). 

Now $\chix_E$ is constant on $K_0\sim K_1$, an open set; since $K_1$ 
is nowhere dense in $K$, given $x\in K_1$, there exists a sequence 
$(x_m)$ in $K_0\sim K_1$ with $x_m\to x$. 
But then $(\osc \chix_E) (x) \ge \lim_{m\to\infty} (\chix_E(x_m) - \chix_E 
(x))=1$, hence (50) is proved for $j=0$. 

Suppose now (50) is proved for $0\le j<n$. 
Again if $x\in K_{j+1}$, since $K_{j+1}$ is nowhere dense in $K_j$, 
choose a sequence $(x_m)$ in $K_j$ with $x_m\to x$. 
Now by definition of $E$, $|\chix_E(x_m)-\chix_E(x)| =1$ for all $m$. 
Thus $\osc \chix_E \mid K_j (x)\ge1$, which proves that $K_{j+1} \subset 
\os_{j+1}(\chix_E,\varep)$. 
But $\chix_E$ is constant on $K_j\sim K_{j+1}$, whence $K_{j+1}\supset 
\osc_{j+1}(\chix_E,\varep)$. 
Thus (50) holds.

To see the final inequality in 2.10, we have that 
$\|\chix_{K_0\sim K_1}\|_D=1$ and \hfill\break 
$\|\chix_{K_{2i}\sim K_{2i+1}}\|_D \le2$ 
for all $1\le i\le [n/2]$ (by Lemma~2.4); hence 
$$\eqalignno{\|\chix_E\|_D &\le  \sum_{i=0}^{[n/2]} \|\chix_{K_{2i}\sim 
K_{2i+1}}\|_D\cr 
&\le 1+ 2 [n/2] \le n+1&\qed\cr}$$ 
\enddemo 

\demo{Remark} 
Actually the final inequality in 2.10 follows  from (49). 
In fact it is proved in \cite{R2} that if $E\subset K$ is such that 
$i(\chix_E)=n$, then $\|\chix_E\|_D=n$ or $n+1$ (and both possibilities 
can occur). 
\enddemo 

\proclaim{Lemma 2.11} 
{\rm (a)} Let $n\ge1$ and suppose $K^{(n)} \ne\emptyset$. There exist 
non-empty closed sets $K_1,\ldots,K_n$ satisfying the hypotheses of 
Lemma~2.10. 
\vskip1pt
{\rm (b)} Suppose $K^{(n)}\ne\emptyset$ for all $n=1,2,\ldots$. 
There exist disjoint open subsets $U_1,U_2,\ldots$ of $K$ with 
$U_n^{(n)} \ne \emptyset$ for all $n$.
\endproclaim 

\demo{Proof} 

(a) If $K$ is perfect, it can be seen that there exists a closed perfect 
nowhere dense subset $L$ of $K$; we then easily obtain the desired sets 
$(K_j)$ with $K_j$ a perfect nowhere dense result of $K_{j-1}$. 
Evidently the same reasoning holds if $K$ has a perfect non-empty subset. 
Otherwise, simply let $K_j = K^{(j)}$, $1\le j\le n$. 
Alternatively, we may just observe that the hypotheses imply $K$ has a closed 
subset homeomorphic to $\omega^n+1$. 

(b) First note that if $x\in K^{(n)}$, then 
$$x\in U^{(n)}\ \text{ for all open neighborhoods }\ U\ \text{ of }\  x\ .
\leqno(51)$$ 
Next, note that the hypotheses imply that $K^{(n)}$ is infinite for all $n$. 
We may thus choose distinct points $x_1,x_2,\ldots,$ with $x_n\in K^{(n)}$ 
for all $n$. 
Now it follows that if $U$ is an open set containing infinitely many of 
the $x_j$'s, there exists an $n$ and an open neighborhood $V$ of $x_n$ with 
$\bar V \subset U$ so that $U\sim\bar V$ contains infinitely many of the 
$x_j$'s. 
We may then  choose $k_1<k_2<\cdots$ and 
$U_1,U_2,\ldots$ open sets with $\bar U_i\cap \bar U_j=\emptyset$ for 
all $i\ne j$ and 
$x_{k_n} \in U_n$ for all $n$. (51) then yields that (b) holds.\qed 
\enddemo 

We finally observe the following simple ``localization'' property for 
$D$-functions. 

\proclaim{Lemma 2.12} 
Let $U_1,U_2,\ldots$ be disjoint non-empty open subsets of $K$, 
$U=\bigcup_{j=1}^\infty U_j$, $\lambda<\infty$, and $f:K\to\real$ a function 
supported on $U$ with $\|f\mid U_j\|_D\le\lambda$ for all $j$. 
Then $f\in DBSC(K)$ and $\|f\|_D \le \lambda$. 
\endproclaim 

\demo{Proof} 
Let $\varep>0$. For each $j$, choose a sequence of continuous functions 
on $K$, $(\varphi_i^j)_{i=1}^\infty$, with $0\le\varphi_i^j \le1$ for all $i$  
and $\chix_{U_j} = \sum_{i=1}^\infty \varphi_i^j$ pointwise. 
Also, choose $(h_i^j)_{i=1}^\infty$ continuous functions on $U_j$, with 
$\sum |h_i^j| \le\lambda+\varep$ and $f\mid U_j = \sum h_i^j$ pointwise. 
Now let 
$$f_{jk\ell} = \varphi_k^j h_\ell^j \chix_{U_j} \ \text{ for all }\ 
j,k,\ell\ .
\leqno(52)$$ 
Then $f_{jk\ell}$ is continuous on $K$ since $h_\ell^j$ is bounded 
continuous on $K$ and supported on $U_j$, and 
$$\displaylines{
\sum_{j,k,\ell} |\varphi_k^j h_\ell^j \chix_{U_j}| 
= \sum_j \sum_\ell |h_\ell^j| \chix_{U_j} \le\lambda +\varep\ ,\cr
\sum_j \sum_\ell \sum_k \varphi_k^j h_\ell^j \chix_{U_j} 
= \sum_j \sum_\ell h_\ell^j \chix_{U_j} 
= \sum_j f\chix_{U_j} = f\ .\cr}$$ 
Thus $\|f\|_D \le \lambda +\varep$; since $\varep>0$ is arbitrary, 
the result follows.\qed
\enddemo 

We are now prepared for the 

\demo{Proof of Proposition 1.5} 

By Lemmas 2.10 and 2.11, we may choose disjoint non-empty open subsets 
$U_1,U_2,\ldots$ of $K$, and for each $n$ a subset $E_n$ of $U_n$ so that 
$$i(\chix_{E_n}) = n = i(\chix_{E_n},\varep)\ \text{ for all }\ 
0< \varep \le1\ .
\leqno(53)$$ 
and 
$$\|\chix_{E_n}\|_{D(U_n)} \le n+1\ .
\leqno(54)$$ 

Now let $f= \sum_{n=1}^\infty \chix_{E_n}/n$ pointwise. 
Thus by Lemma~2.12 and (54), $f\in DBSC(K)$ (with $\|f\|_D\le 2$). 
However fixing $n$ and letting $\varep=\tfrac1n$, then by (53), 
$i(\chix_{E_n},1) =n$ ($= i({1\over n} \chix_{E_n}, {1\over n})$) and so 
$$\varep i(f,\varep) \ge \frac1n i\Bigl( f\mid U_n,\frac1n\Bigr) =1\ .
\leqno(55)$$ 
Thus $f$ fails (7), so $f\notin SD(K)$ by 
Theorem 1.4.\qed 
\enddemo

\bigskip
\baselineskip=12pt 
\leftline{\sl Authors addresses:} 
\medskip

{\vbox{\halign{{#}\hfill\qquad\qquad&{#}\hfill\cr
F. Chaatit&V. Mascioni and H. Rosenthal\cr
Department of Mathematics&Department of Mathematics\cr
University of Texas at El Paso&University of Texas at Austin\cr
El Paso, TX 79968-0514&Austin, TX 78712-1082\cr}}}

\enddocument